\documentclass[11pt,a4paper]{amsart}

\usepackage[T1]{fontenc}
\usepackage{lmodern}
\usepackage{amsmath,amssymb,amsthm}
\usepackage[hidelinks]{hyperref}

\newtheorem{theorem}{Theorem}[section]
\newtheorem{proposition}[theorem]{Proposition}
\newtheorem{lemma}[theorem]{Lemma}
\newtheorem{corollary}[theorem]{Corollary}
\theoremstyle{plain}
\newtheorem*{conjecture}{Conjecture}
\numberwithin{equation}{section}

\newcommand{\N}{\mathbb N}
\newcommand{\Pp}{\mathbb P}
\newcommand{\Pc}{\mathcal P}
\newcommand{\Nc}{\mathcal N}
\newcommand{\Dc}{\mathcal D}

\title[Cram\'er's conjecture and Beurling's conjecture]{Cram\'er's Conjecture and Beurling's Conjecture Are Incompatible}
\author{Volodymyr Pavlenkov}
\address{Department of Mathematics, University of York,
Heslington, York, YO10 5DD, England, {\texttt{volodymr.pavlenkov@york.ac.uk}}}

\begin{document}

\begin{abstract}
Beurling conjectured that a system of generalized primes with associated integer-counting function $N(x)=x+o(\log x)$ must coincide with the classical prime system. We construct a nonclassical system satisfying this asymptotic relation under the assumption that the gaps between consecutive ordinary primes are $O(p_n^\theta)$ for some $\theta<1/2$. In particular, Cram\'er's conjecture implies that Beurling's conjecture is false. Conversely, the validity of Beurling's conjecture would imply that, for every $\varepsilon>0$, the ratios $(p_{n+1}-p_n)/p_n^{1/2-\varepsilon}$ are unbounded.
\end{abstract}

\maketitle

\section{Introduction}

Let $\Pc=\{q_1,q_2,\ldots\}$ be a sequence of real numbers satisfying
\[
  1<q_1<q_2<\cdots,\qquad q_j\to\infty.
\]
The elements of $\Pc$ are called Beurling generalized primes. The associated set of generalized integers is defined by
\[
  \Nc_{\Pc}=\left\{\prod_{j\geq 1}q_j^{\nu_j}:
  \nu_j\in\mathbb Z_{\geq0},\ \nu_j=0\text{ for all but finitely many }j\right\},
\]
and we denote its counting function by
\[
  N_{\Pc}(x)=\#\{n\in\Nc_{\Pc}:n\leq x\}.
\]
Let $\Pp=\{p_1,p_2,\ldots\}$ be the set of ordinary primes, enumerated in increasing order. By the fundamental theorem of arithmetic, $\Nc_{\Pp}=\N$ and $N_{\Pp}(x)=\lfloor x\rfloor$.

The following conjecture of Beurling is stated in Erd\H{o}s' survey on combinatorial number theory \cite[p.~103, (12)]{Erdos1980}; see also \cite[Problem~951]{Bloom}.

\begin{conjecture}[Beurling]
Suppose that distinct formal products of elements of $\Pc$ have distinct numerical values. If
\begin{equation}\label{eq:beurling}
  N_{\Pc}(x)=x+o(\log x),\qquad x\to\infty,
\end{equation}
then $\Pc=\Pp$.
\end{conjecture}

In \cite{Erdos1980} the counting function is defined using the strict inequality $n<x$. For a set of generalized integers the two counting conventions differ by at most one, and thus give the same conjecture. The uniqueness assumption makes explicit the convention needed when counting distinct numerical values: without it, adjoining a composite integer to the ordinary primes would leave the counting function unchanged. No separation condition on the generalized integers is imposed. In the system constructed below all distinct formal products have distinct numerical values, so counting with multiplicities does not change its counting function.

Another conjecture concerning the ordinary primes, due to Cram\'er \cite{Cramer1936}, is commonly expressed in the following upper-bound form.

\begin{conjecture}[Cram\'er]
The gaps between consecutive primes satisfy
\begin{equation}\label{eq:cramer}
  p_{n+1}-p_n=O((\log p_n)^2),\qquad n\to\infty.
\end{equation}
\end{conjecture}

We shall show that the two conjectures cannot hold simultaneously. More precisely, we establish the following result under a substantially weaker assumption on the gaps between consecutive ordinary primes.

\begin{theorem}\label{thm:main}
Suppose that for some $0<\theta<1/2$,
\begin{equation}\label{eq:gaps}
  p_{n+1}-p_n=O(p_n^\theta),\qquad n\to\infty.
\end{equation}
Then there exists a generalized prime system $\Pc^*\ne\Pp$ such that
\begin{equation}\label{eq:main}
  N_{\Pc^*}(x)=x+O\big((\log x)^{\beta(\theta)}\big),\qquad x\to\infty,
\end{equation}
where
\[
  \beta(\theta)=\frac{\log 2}{\log(3-2\theta)}<1.
\]
In particular, Beurling's conjecture is false under \eqref{eq:gaps}.
\end{theorem}

The system in Theorem~\ref{thm:main} is obtained by deleting the primes $2,3,5$ and replacing a sufficiently sparse sequence of ordinary primes by their square roots. In Section~2 we derive an estimate for the integer-counting function of such a system. In Section~3 we construct the required sequence of ordinary primes under \eqref{eq:gaps} and deduce the consequences for the two conjectures. We also obtain estimates from the unconditional prime-gap bound of Baker, Harman and Pintz and from the classical bound under the Riemann hypothesis.

\section{A construction of generalized primes}

Let $5<r_1<r_2<\cdots$ be a sequence of ordinary primes and consider the system
\begin{equation}\label{eq:system}
  \Pc^*=\big(\Pp\setminus\{2,3,5,r_1,r_2,\ldots\}\big)
         \cup\{\sqrt{r_1},\sqrt{r_2},\ldots\}.
\end{equation}
Denote by $\Dc$ the set of square-free integers of the form
\[
  d=\prod_{j\in J}r_j,\qquad J\subset\N,\quad |J|<\infty,
\]
where the empty product is equal to one. Define
\[
  B(y)=\#\{m\in\N:m\leq y,\ (m,30)=1\},\qquad y\geq0.
\]
By the inclusion--exclusion principle,
\begin{equation}\label{eq:B}
  B(y)=\sum_{a\mid30}\mu(a)\left\lfloor\frac{y}{a}\right\rfloor
      =\frac{4}{15}y+O(1),
\end{equation}
where the error term is uniform for $y\geq0$ and $\mu$ denotes the M\"obius function.

\begin{lemma}\label{lem:count}
For the system defined by \eqref{eq:system},
\begin{equation}\label{eq:count}
  N_{\Pc^*}(x)=\sum_{\substack{d\in\Dc\\d\leq x^2}}
                B\left(\frac{x}{\sqrt d}\right).
\end{equation}
Moreover, if $\sqrt{r_k}\leq x<\sqrt{r_{k+1}}$, then
\begin{equation}\label{eq:countfinite}
  N_{\Pc^*}(x)=\frac{4}{15}x\prod_{j=1}^k(1+r_j^{-1/2})+O(2^k),
\end{equation}
where the implied constant is independent of $x$ and $k$.
\end{lemma}

\begin{proof}
Every non-negative integral power of $\sqrt{r_j}$ is of the form $r_j^a(\sqrt{r_j})^b$, where $a\geq0$ is an integer and $b\in\{0,1\}$. Since $r_j=(\sqrt{r_j})^2$, all ordinary primes other than $2,3,5$ belong to $\Nc_{\Pc^*}$. Consequently, the generalized integers are precisely the numbers
\begin{equation}\label{eq:representation}
  m\sqrt d,\qquad m\in\N,\quad (m,30)=1,\quad d\in\Dc.
\end{equation}
If $m\sqrt d=m'\sqrt{d'}$, then $m^2d=(m')^2d'$. Comparing the parities of the prime exponents on both sides, we obtain $d=d'$ and hence $m=m'$. The fundamental theorem of arithmetic also implies that each number in \eqref{eq:representation} has a unique representation as a formal product of the elements of $\Pc^*$: for every $r_j$ the exponent of $\sqrt{r_j}$ is $2v_{r_j}(m)+\mathbf 1_{r_j\mid d}$, where $v_p(m)$ is the exponent of $p$ in $m$. This proves \eqref{eq:count}.

Suppose that $\sqrt{r_k}\leq x<\sqrt{r_{k+1}}$. Any $d\in\Dc$ divisible by some $r_j$ with $j>k$ satisfies $\sqrt d>x$. Since $B(y)=0$ for $y<1$, the sum in \eqref{eq:count} can therefore be written as
\[
  N_{\Pc^*}(x)=\sum_{J\subset\{1,\ldots,k\}}
     B\left(x\prod_{j\in J}r_j^{-1/2}\right).
\]
Applying \eqref{eq:B} to all $2^k$ terms, we obtain
\begin{align*}
  N_{\Pc^*}(x)
    &=\frac{4}{15}x\sum_{J\subset\{1,\ldots,k\}}
       \prod_{j\in J}r_j^{-1/2}+O(2^k)\\
    &=\frac{4}{15}x\prod_{j=1}^k(1+r_j^{-1/2})+O(2^k).
\end{align*}
\end{proof}

\begin{proposition}\label{prop:sparse}
Suppose that the sequence $(r_j)$ satisfies
\begin{equation}\label{eq:normalization}
  \prod_{j=1}^{\infty}(1+r_j^{-1/2})=\frac{15}{4}
\end{equation}
and that for some $c>0$ and $\lambda>1$,
\begin{equation}\label{eq:sparse}
  r_{j+1}\geq c r_j^{\lambda}
\end{equation}
for all sufficiently large $j$. Then the system defined by \eqref{eq:system} satisfies
\begin{equation}\label{eq:error}
  N_{\Pc^*}(x)=x+O\big((\log x)^{\log 2/\log\lambda}\big),
  \qquad x\to\infty.
\end{equation}
Furthermore, if $\lambda>2$, then
\begin{equation}\label{eq:L}
  \liminf_{x\to\infty}(\pi_{\Pc^*}(x)-\pi(x))=-3.
\end{equation}
where $\pi_{\Pc^*}$ and $\pi$ denote the generalized and ordinary prime-counting functions, respectively.
\end{proposition}

\begin{proof}
By \eqref{eq:normalization},
\[
  \frac{4}{15}\prod_{j=1}^k(1+r_j^{-1/2})
   =\exp(-T_k),\qquad
  T_k=\sum_{j>k}\log(1+r_j^{-1/2}).
\]
For $k$ sufficiently large, condition \eqref{eq:sparse} implies that $r_{j+1}\geq4r_j$ for every $j\geq k+1$. Thus
\[
  0\leq T_k\leq\sum_{j>k}r_j^{-1/2}\leq 2r_{k+1}^{-1/2}.
\]
Taking $k$ such that $\sqrt{r_k}\leq x<\sqrt{r_{k+1}}$, we have $xT_k\leq2$. Since $0\leq1-e^{-T_k}\leq T_k$, Lemma~\ref{lem:count} implies
\begin{equation}\label{eq:error2k}
  |N_{\Pc^*}(x)-x|\ll 2^k+xT_k\ll 2^k.
\end{equation}

Set $u_j=\log r_j$ and $C=\max\{0,-\log c\}$. Taking logarithms in \eqref{eq:sparse}, we obtain $u_{j+1}\geq\lambda u_j-C$ for all sufficiently large $j$. Choose $j_0$ such that this inequality holds for all $j\geq j_0$ and $u_{j_0}>C/(\lambda-1)$. Iterating the inequality, we obtain
\[
  u_k-\frac{C}{\lambda-1}\geq
    \lambda^{k-j_0}\left(u_{j_0}-\frac{C}{\lambda-1}\right),\qquad k\geq j_0.
\]
In particular, there exists $c_1>0$ such that $\log r_k\geq c_1\lambda^k$ for every sufficiently large $k$. Since $x\geq\sqrt{r_k}$,
\[
  2^k=(\lambda^k)^{\log 2/\log\lambda}
     \ll (\log x)^{\log 2/\log\lambda}.
\]
Combining this estimate with \eqref{eq:error2k} gives \eqref{eq:error}.

Finally, suppose that $\lambda>2$ and let $A(x)=\#\{j:r_j\leq x\}$. For $x\geq5$, \eqref{eq:system} implies that
\[
  \pi_{\Pc^*}(x)-\pi(x)=-3+A(x^2)-A(x)\geq-3.
\]
Since $\lambda>2$, condition \eqref{eq:sparse} implies that $r_{k+1}>r_k^2$ for every sufficiently large $k$. Hence $A(r_k^2)=A(r_k)=k$ for such $k$, and \eqref{eq:L} follows.
\end{proof}

\section{Gaps between consecutive ordinary primes}

We next construct a sequence of primes satisfying the assumptions of Proposition~\ref{prop:sparse}.

\begin{lemma}\label{lem:greedy}
Suppose that \eqref{eq:gaps} holds for some $0<\theta<1$. Then, for every $a>0$, there exists a strictly increasing sequence $(r_j)$ of ordinary primes greater than $5$ such that
\begin{equation}\label{eq:sum}
  \sum_{j=1}^{\infty}\log(1+r_j^{-1/2})=a
\end{equation}
and, for some $c>0$,
\begin{equation}\label{eq:growth}
  r_{j+1}\geq c r_j^{3-2\theta}
\end{equation}
for every sufficiently large $j$.
\end{lemma}

\begin{proof}
Let $f(t)=\log(1+t^{-1/2})$ for $t>0$. The function $f$ is positive and strictly decreasing, and
\begin{equation}\label{eq:f}
  f(t)\sim t^{-1/2},\qquad
  f'(t)=-\frac{1}{2t(\sqrt t+1)},\qquad t\to\infty.
\end{equation}
Its inverse is given by $f^{-1}(u)=(e^u-1)^{-2}$ for $u>0$.

For $t>1$, let $P(t)$ be the least ordinary prime strictly greater than $t$. By \eqref{eq:gaps}, there exists $C_0>0$ such that
\begin{equation}\label{eq:nextprime}
  0<P(t)-t\leq C_0 t^\theta
\end{equation}
for all sufficiently large $t$. In particular, $P(t)\sim t$. Applying the mean value theorem and \eqref{eq:f}, we obtain
\begin{equation}\label{eq:one-step}
  0<f(t)-f(P(t))\leq \frac{C_0}{2}t^{\theta-3/2}
\end{equation}
for all sufficiently large $t$. Since $\theta<1$, there exists $T>5$ such that, for all $t\geq T$,
\begin{equation}\label{eq:half}
  f(t)-f(P(t))\leq\frac12 f(t),\qquad P(t)\leq 2t,
\end{equation}
and \eqref{eq:one-step} holds.

We first choose a finite sequence of primes. Let $s_1<s_2<\cdots$ be the sequence of all ordinary primes exceeding $T$. Since $f(p)\geq1/(2\sqrt p)$ for $p>1$, the divergence of $\sum_{p\in\Pp}p^{-1}$ implies that $\sum_j f(s_j)=\infty$. Thus there exists a smallest integer $M\geq1$ such that
\[
  \sum_{j=1}^{M}f(s_j)\geq a.
\]
Put $m=M-1$, set $r_j=s_j$ for $1\leq j\leq m$, and let
\[
  \delta_m=a-\sum_{j=1}^m f(r_j).
\]
Here the sum is empty if $m=0$. By the choice of $M$, we have $0<\delta_m\leq f(s_M)$. In particular,
\[
  t_m:=f^{-1}(\delta_m)\geq s_M>T,
\]
and $t_m>r_m$ if $m\geq1$.

Suppose that $r_1,\ldots,r_k$ have been chosen, where $k\geq m$, and that
\[
  \delta_k=a-\sum_{j=1}^k f(r_j)>0,\qquad t_k=f^{-1}(\delta_k)\geq T.
\]
Define recursively
\begin{equation}\label{eq:recursive}
  r_{k+1}=P(t_k),\qquad
  \delta_{k+1}=\delta_k-f(r_{k+1})=f(t_k)-f(P(t_k)).
\end{equation}
By \eqref{eq:half}, $0<\delta_{k+1}\leq\delta_k/2$. In particular,
\[
  \delta_{k+1}\leq\delta_k-\delta_{k+1}=f(r_{k+1}).
\]
Consequently, $t_{k+1}\geq r_{k+1}>t_k$, and the next choice $r_{k+2}=P(t_{k+1})$ satisfies $r_{k+2}>r_{k+1}$. The recursive construction is therefore well defined and gives a strictly increasing sequence of primes. Also, $\delta_k\to0$, which proves \eqref{eq:sum}.

By \eqref{eq:one-step} and $r_{k+1}\leq2t_k$, there exists $C_1>0$ such that
\[
  \delta_{k+1}\leq C_1 r_{k+1}^{\theta-3/2}.
\]
Since $f^{-1}(u)\sim u^{-2}$ as $u\downarrow0$, for all sufficiently large $k$ we have
\begin{align*}
  r_{k+2}>t_{k+1}=f^{-1}(\delta_{k+1})
     &\geq\frac{1}{2}\delta_{k+1}^{-2}\\
     &\geq\frac{1}{2C_1^2} r_{k+1}^{3-2\theta}.
\end{align*}
This proves \eqref{eq:growth}.
\end{proof}

\begin{proof}[Proof of Theorem~\ref{thm:main}]
Apply Lemma~\ref{lem:greedy} with $a=\log(15/4)$. The resulting sequence satisfies \eqref{eq:normalization} and \eqref{eq:sparse} with $\lambda=3-2\theta>2$. By Proposition~\ref{prop:sparse}, the system \eqref{eq:system} satisfies \eqref{eq:main}. Since the prime $2$ is absent from $\Pc^*$, we have $\Pc^*\ne\Pp$. Finally, $\beta(\theta)<1$, and hence \eqref{eq:beurling} holds for $\Pc^*$.
\end{proof}

\begin{corollary}\label{cor:cramer}
Cram\'er's conjecture and Beurling's conjecture are incompatible. More precisely, under Cram\'er's conjecture, for every $\beta>\log 2/\log 3$ there exists a generalized prime system $\Pc^*\ne\Pp$ such that
\[
  N_{\Pc^*}(x)=x+O((\log x)^\beta),\qquad x\to\infty.
\]
\end{corollary}

\begin{proof}
Under \eqref{eq:cramer}, condition \eqref{eq:gaps} holds for every $\theta>0$. Since $\beta(\theta)\to\log 2/\log 3$ as $\theta\downarrow0$, for each $\beta>\log 2/\log 3$ we may choose $\theta\in(0,1/2)$ such that $\beta(\theta)<\min\{1,\beta\}$. The assertion follows from Theorem~\ref{thm:main}.
\end{proof}

The contrapositive of Theorem~\ref{thm:main} gives the following necessary condition for Beurling's conjecture to hold.

\begin{corollary}\label{cor:gaps}
If Beurling's conjecture is true, then for every $\varepsilon>0$,
\begin{equation}\label{eq:largegaps}
  \limsup_{n\to\infty}\frac{p_{n+1}-p_n}{p_n^{1/2-\varepsilon}}=+\infty.
\end{equation}
\end{corollary}

\begin{proof}
Suppose that \eqref{eq:largegaps} fails for some $\varepsilon>0$. Then
\[
  p_{n+1}-p_n=O(p_n^{1/2-\varepsilon}).
\]
Choose $\theta$ such that $\max\{0,1/2-\varepsilon\}<\theta<1/2$. The hypothesis \eqref{eq:gaps} of Theorem~\ref{thm:main} is satisfied, and hence there exists a nonclassical generalized prime system satisfying \eqref{eq:beurling}. This contradicts Beurling's conjecture.
\end{proof}

We finally record the estimates obtained from the same construction using the unconditional prime-gap bound of Baker, Harman and Pintz and the classical bound under the Riemann hypothesis.

\begin{corollary}\label{cor:bhp}
There exists a generalized prime system $\Pc^*\ne\Pp$ of the form \eqref{eq:system} such that
\[
  N_{\Pc^*}(x)=x+O\big((\log x)^{\log 2/\log(39/20)}\big),\qquad x\to\infty.
\]
Here $\log 2/\log(39/20)=1.037910\ldots$.
\end{corollary}

\begin{proof}
By the theorem of Baker, Harman and Pintz \cite{BHP2001},
\[
  p_{n+1}-p_n\ll p_n^{21/40}.
\]
Apply Lemma~\ref{lem:greedy} with $a=\log(15/4)$ and $\theta=21/40$. The resulting sequence satisfies \eqref{eq:normalization} and \eqref{eq:sparse} with $\lambda=39/20$. The assertion follows from Proposition~\ref{prop:sparse}.
\end{proof}

\begin{corollary}\label{cor:rh}
Assume the Riemann hypothesis. Then there exists a generalized prime system $\Pc^*\ne\Pp$ of the form \eqref{eq:system} such that
\[
  N_{\Pc^*}(x)=x+O(\log x),\qquad x\to\infty.
\]
\end{corollary}

\begin{proof}
Under the Riemann hypothesis, the classical estimate of Cram\'er gives
\[
  p_{n+1}-p_n\ll \sqrt{p_n}\log p_n.
\]
Consequently, for the least prime $P(t)$ strictly greater than $t$, we have $P(t)-t\ll \sqrt t\log t$. By the mean value theorem and \eqref{eq:f},
\[
  0<f(t)-f(P(t))\ll t^{-1}\log t.
\]
In particular, \eqref{eq:half} holds for all sufficiently large $t$. Repeating the construction in the proof of Lemma~\ref{lem:greedy} with $a=\log(15/4)$, we obtain \eqref{eq:normalization} and
\[
  \delta_{k+1}\ll \frac{\log r_{k+1}}{r_{k+1}}.
\]
Since $f^{-1}(u)\sim u^{-2}$ as $u\downarrow0$, there exists $c>0$ such that, for all sufficiently large $k$,
\begin{equation}\label{eq:rhgrowth}
  r_{k+1}\geq c\frac{r_k^2}{(\log r_k)^2}.
\end{equation}
Set $u_k=\log r_k$. By \eqref{eq:rhgrowth}, for some $C\geq0$,
\[
  u_{k+1}\geq 2u_k-2\log u_k-C
\]
for all sufficiently large $k$. Put $D=C+2\log 2$ and $g(u)=u-2\log u-D$. The function $g$ is increasing for $u>2$. For all sufficiently large $k$, the number $v_k=2u_k-2\log u_k-C$ satisfies $2<v_k\leq2u_k$. Thus
\begin{align*}
  g(u_{k+1})&\geq g(v_k)\\
  &\geq 2u_k-2\log u_k-C-2\log(2u_k)-D\\
  &=2g(u_k).
\end{align*}
Choose $k_0$ sufficiently large that this inequality holds for all $k\geq k_0$ and $g(u_{k_0})>0$. Iteration gives
\[
  u_k\geq g(u_k)\geq 2^{k-k_0}g(u_{k_0})\gg 2^k.
\]
Moreover, \eqref{eq:rhgrowth} implies that $r_{k+1}\geq4r_k$ for all sufficiently large $k$. The argument leading to \eqref{eq:error2k} therefore applies. For $\sqrt{r_k}\leq x<\sqrt{r_{k+1}}$, we obtain
\[
  |N_{\Pc^*}(x)-x|\ll 2^k\ll\log r_k\leq 2\log x.
\]
This proves the assertion.
\end{proof}

Neither of the last two estimates gives the little-$o$ relation \eqref{eq:beurling}. In particular, the bound obtained from the Riemann hypothesis reaches the logarithmic order but does not by itself disprove Beurling's conjecture.


\begin{thebibliography}{9}

\bibitem{BHP2001}
R.~C. Baker, G.~Harman and J.~Pintz, The difference between consecutive primes, II, \emph{Proc. London Math. Soc.} (3) \textbf{83} (2001), 532--562.

\bibitem{Bloom}
T.~F. Bloom, \emph{Erd\H{o}s Problems}, Problem~951,
\url{https://www.erdosproblems.com/951}.

\bibitem{Cramer1936}
H.~Cram\'er, On the order of magnitude of the difference between consecutive prime numbers, \emph{Acta Arith.} \textbf{2} (1936), 23--46.

\bibitem{Erdos1980}
P.~Erd\H{o}s, A survey of problems in combinatorial number theory, \emph{Ann. Discrete Math.} \textbf{6} (1980), 89--115.

\end{thebibliography}
\end{document}